\documentclass[11pt,reqno]{amsart}
\usepackage{amssymb,verbatim,enumerate,ifthen}
\usepackage[mathscr]{eucal}
\usepackage[utf8]{inputenc}
\usepackage[T1]{fontenc}
\allowdisplaybreaks[1]

\def\R{\mathbb{R}}

\newtheorem{theorem}{Theorem}
\newtheorem*{theorem*}{Theorem}
\def\Thm#1#2{\ifthenelse{\equal{#1}{*}}{\begin{theorem*}#2\end{theorem*}}
             {\begin{theorem}\label{T#1}#2\end{theorem}}}
\newtheorem{Atheorem}{Theorem}

\def\thm#1{Theorem~\ref{T#1}}

\newtheorem{proposition}[theorem]{Proposition}
\newtheorem*{proposition*}{Proposition}
\def\Prp#1#2{\ifthenelse{\equal{#1}{*}}{\begin{proposition*}#2\end{proposition*}}
             {\begin{proposition}\label{P#1}#2\end{proposition}}}

\newtheorem{corollary}[theorem]{Corollary}
\newtheorem*{corollary*}{Corollary}
\def\Cor#1#2{\ifthenelse{\equal{#1}{*}}{\begin{corollary*}#2\end{corollary*}}
             {\begin{corollary}\label{C#1}#2\end{corollary}}}
\def\cor#1{Corollary~\ref{C#1}}

\newtheorem{lemma}[theorem]{Lemma}
\newtheorem*{lemma*}{Lemma}
\def\Lem#1#2{\ifthenelse{\equal{#1}{*}}{\begin{lemma*}#2\end{lemma*}}
             {\begin{lemma}\label{L#1}#2\end{lemma}}}
\def\lem#1{Lemma~\ref{L#1}}
\newtheorem{Alemma}{Lemma}

\theoremstyle{definition}
\newtheorem{remark}[theorem]{Remark}
\newtheorem*{remark*}{Remark}
\def\Rem#1#2{\ifthenelse{\equal{#1}{*}}{\begin{remark*}\rm #2\end{remark*}}
             {\begin{remark}\label{R#1}\rm #2\end{remark}}}

\newtheorem{example}[theorem]{Example}
\newtheorem*{example*}{Example}
\def\Exa#1#2{\ifthenelse{\equal{#1}{*}}{\begin{example*}\rm #2\end{example*}}
             {\begin{example}\label{Ex#1}\rm #2\end{example}}}

\def\eq#1{{\rm(\ref{E#1})}}
\def\Eq#1#2{\ifthenelse{\equal{#1}{*}}
  {\begin{equation*}\begin{aligned}#2\end{aligned}\end{equation*}}
  {\begin{equation}\begin{aligned}\label{E#1}#2\end{aligned}\end{equation}}}

\def\diag{\mathop{\hbox{\rm diag}}\nolimits}

\begin{document}
\vspace{5mm}

\date{\today}

\title[Comparison of nonsymmetric generalized Bajraktarevi\'c means]{Local and global comparison of nonsymmetric generalized Bajraktarevi\'c means}

\author[R.\ Gr\"unwald]{Rich\'ard Gr\"unwald}
\address[R.\ Gr\"unwald]{Doctoral School of Mathematical and Computational Sciences, University of Debrecen, H-4002 Debrecen, Pf. 400, Hungary}
\email{richard.grunwald@science.unideb.hu}

\author[Zs. P\'ales]{Zsolt P\'ales}
\address[Zs. P\'ales]{Institute of Mathematics, University of Debrecen, H-4002 Debrecen, Pf. 400, Hungary}
\email{pales@science.unideb.hu}

\thanks{The research of the first author was supported by the \'UNKP-21-3 New National Excellence Program of the Ministry of Human Capacities. The research of the second author was supported by the K-134191 NKFIH Grant and the 2019-2.1.11-TÉT-2019-00049, the EFOP-3.6.1-16-2016-00022 project. The last project is co-financed by the European Union and the European Social Fund.}


\begin{abstract}
The purpose of this paper is to investigate the local and global comparison of two $n$-variable nonsymmetric generalized Bajraktarević means, i.e., to establish necessary as well as sufficient conditions in terms of the unknown functions $f,g,p_1,\dots,p_n,q_1,\dots,q_n:I\to\mathbb{R}$ for the comparison inequality
$$
f^{-1}\bigg(\frac{p_1(x_1)f(x_1)+\cdots+p_n(x_n)f(x_n)}{p_1(x_1)+\cdots+p_n(x_n)}\bigg)\leq g^{-1}\bigg(\frac{q_1(x_1)g(x_1)+\cdots+q_n(x_n)g(x_n)}{q_1(x_1)+\cdots +q_n(x_n)}\bigg)
$$
in local and global sense. Here $I$ is a nonempty open real interval, $x_1,\dots,x_n\in I$, and $f,g$ are assumed to be continuous, strictly monotone and $p_1,\dots,p_n,q_1,\dots,q_n:I\to\mathbb{R}_+$ are positive valued. Concerning the global comparison problem, the main result of the paper states that if $f,g$ are differentiable functions with nonvanishing first derivatives and, for all $i\in\{1,\dots,n\}$,
$$
 \frac{p_i}{p_0}=\frac{q_i}{q_0} \qquad\mbox{and}\qquad
\frac{p_0(x)(f(x)-f(y))}{p_0(y)f'(y)}
\leq\frac{q_0(x)(g(x)-g(y))}{q_0(y)g'(y)}\qquad(x,y\in I)
$$
are satisfied (where $p_0:=p_1+\dots+p_n$ and $q_0:=q_1+\dots+q_n$), then the above comparison inequality holds for all $x_1,\dots,x_n\in I$.
\end{abstract}

\maketitle

\section{Introduction}

Throughout this paper, the symbols $\R$ and $\R_+$ will stand for the sets of real and positive real numbers, respectively, and $I$ will always denote a nonempty open real interval.

Given a strictly monotone continuous function $f:I\to\R$ and an $n$-tuple of positive valued functions $p=(p_1,\dots,p_n):I\to\R_{+}^n$, the \emph{$n$-variable nonsymmetric generalized Bajraktarević mean} $A_{f,p}:I^n\to I$ is given by the following formula:
\Eq{BM}{
	A_{f,p}(x_1,\dots,x_n):=f^{(-1)}\bigg(\frac{p_1(x_1)f(x_1)+\dots+p_n(x_n)f(x_n)}{p_1(x_1)+\dots+p_n(x_n)}\bigg) \qquad (x_1,\dots,x_n\in I).
}
This is an extension of the notion introduced by Bajraktarevi\'c in the symmetric setting in \cite{Baj58} and \cite{Baj63} which was the case when $p_1=\dots=p_n$, i.e., when all the weight functions are the same. In the sequel, the sum of these weight functions will be denoted by $p_0$, i.e., $p_0:=p_1+\dots+p_n$. It is easy to see that $A_{f,p}$ is a strict mean, i.e.,
\Eq{*}{
	\min(x_1,\dots,x_n)\leq A_{f,p}(x_1,\dots,x_n)\leq \max(x_1,\dots,x_n)\qquad (x_1,\dots,x_n\in I)
}
holds and the inequalities are strict if $\min(x_1,\dots,x_n)<\max(x_1,\dots,x_n)$.

The equality problem of nonsymmetric generalized Bajraktarević means has been solved by the authors in \cite{GruPal20}. The main goal of this paper is to investigate the local and global comparison problem of two $n$-variable nonsymmetric generalized Bajraktarević means, i.e., to find necessary as well as sufficient conditions for the functional inequality
\Eq{comp}{
	A_{f,p}(x_1,\dots,x_n)\leq A_{g,q}(x_1,\dots,x_n),
}
where $f,g:I\to\R$ are strictly monotone functions and $p,q:I\to\R_+^n$. If there exists an open set $U\subseteq I^n$ which contains the diagonal of $I^n$ and \eq{comp} holds for all $(x_1,\dots,x_n)\in U$, then we say that \emph{$A_{f,p}$ is locally smaller than $A_{g,q}$}. If \eq{comp} holds for all $(x_1,\dots,x_n)\in I^n$, then we say that \emph{$A_{f,p}$ is globally smaller than $A_{g,q}$}. Clearly, the global comparability of the means implies their local comparability.

The (global) comparison problem of classical Bajraktarević means (with arbitrary number of variables) was solved by Daróczy and Losonczi \cite{DarLos70}. Their result was extended to the context of deviation means by Daróczy \cite{Dar72b} and by Daróczy and Páles \cite{DarPal82}. The so-called Gini means (cf.\ \cite{Gin38}) according to the result of Aczél and Daróczy \cite{AczDar63c} are precisely the continuous and homogeneous Bajraktarević means on $I=\R_+$. Their comparison was also characterized by Daróczy and Losonczi \cite{DarLos70}. For the setting when $I$ is a proper subinterval of $\R_+$, the comparison problem of Gini means was solved by Losonczi \cite{Los77} and this result was extended to the setting of Gini means with complex parameters by Páles \cite{Pal89c}. Another generalization of Bajraktarević means (in terms of a measure and two functions) was introduced and their comparison problem was investigated by Losonczi and Páles \cite{LosPal08}. The notion of local and global comparison in broad classes of means was introduced by Páles and Zakaria \cite{PalZak17}, where such comparison problems were in the focus of research. Motivated by theses preliminaries, in this paper we establish necessary as well as sufficient conditions both for the local and for the global comparability of nonsymmetric generalized Bajraktarević means. We note that, for nonsymmetric generalized Bajraktarević means, the equality in local and in global sense are equivalent properties. However, as we shall see, this is not the case for the comparison problem.

\section{Local comparison of Bajraktarević means}

In order to investigate inequality \eq{comp} in the local sense, we recall the following result from \cite{PalZak17}. In the formulation of the statements below, for an $n$ variable function $\Phi:I^n\to\R$, the function $\Phi^\Delta:I\to\R$ is defined by
\Eq{*}{
  \Phi^\Delta(x):=\Phi(x,\dots,x)\qquad(x\in I).
}

\Thm{PZ17}{Let $n\geq 2$ and let $M,N:I^n\to I$ be $n$-variable means such that $M$ is locally smaller than $N$. Assume that $M$ and $N$ are partially differentiable at each diagonal points of $I^n$. Then, for all $i\in\{1,\dots,n\}$,
\Eq{1A}{
  \partial_i M^\Delta=\partial_i N^\Delta.
}
If, in addition, $M$ and $N$ are twice differentiable at each diagonal point of $I^n$, then the symmetric $(n-1)\times (n-1)$-matrix
\Eq{1B}{
  \big(\partial_i\partial_j N^\Delta-\partial_i\partial_j M^\Delta\big)_{i,j=1}^{n-1}
}
is positive semidefinite. \\\indent On the other hand, if the equality \eq{1A} holds for all $i\in\{1,\dots,n\}$, furthermore, $M$ and $N$ are twice continuously differentiable at each point of the diagonal of $I^n$ and the symmetric $(n-1)\times (n-1)$-matrix given by \eq{1B} is positive definite on $I$, then $M$ is locally smaller than $N$.}

In order to apply this result for the nonsymmetric generalized Bajraktarević means, we need to compute their partial derivatives at the diagonal points of $I^n$. For this aim, we recall the following result, which was obtained by the authors in \cite{GruPal20}.

\Lem{DB}{
	Let $\ell\in\{1,2\}$, let $f: I\to\R$ be an $\ell$ times differentiable function on $I$ with a nonvanishing first derivative, and let $p=(p_1,\dots,p_n):I\to\R_{+}^n$. Then we have the following assertions.
	\begin{enumerate}
		\item[(1)] If $\ell=1$, $i\in \{1,\dots,n\}$, and $p_i$ is continuous on $I$, then the first-order partial derivative $\partial_i A_{f,p}$ exists on $\diag(I^n)$ and
		\Eq{*}{
		\partial_i A_{f,p}^\Delta=\frac{p_i}{p_0}.
		}
		\item[(2a)] If $\ell=2$, $i\in \{1,\dots,n\}$, and $p_i$ is continuously differentiable on $I$, then the second-order partial derivative $\partial_i^2 A_{f,p}$ exists on $\diag(I^n)$ and
		\Eq{*}{
		\partial_i^2 A_{f,p}^\Delta=2\frac{p_i'(p_0-p_i)}{p_0^2}+\frac{p_i(p_0-p_i)}{p_0^2}\cdot\frac{f''}{f'}.
		}
		\item[(2b)] If $\ell=2$, $i,j\in \{1,\dots,n\}$ with $i\neq j$, furthermore, $p_i$ and $p_j$ are differentiable on $I$, then the second-order partial derivative $\partial_i\partial_j A_{f,p}$ exists on $\diag(I^n)$ and
		\Eq{*}{
		\partial_i\partial_j A_{f,p}^\Delta =-\frac{(p_ip_j)'}{p_0^2}-\frac{p_ip_j}{p_0^2}\cdot\frac{f''}{f'}.
		}
	\end{enumerate}
	In addition, if $f$ and $p$ are $\ell$ times continuously differentiable, then $A_{f,p}$ is also $\ell$ times continuously differentiable on $I^n$.
}

Our first main result establishes necessary and sufficient conditions for the local comparability of nonsymmetric generalized Bajraktarević means.

\Thm{1NC}{Let $n\geq 2$ and let $f,g:I\to\R$ be differentiable functions with nonvanishing first derivatives and $p,q:I\to\R_+^n$ be continuous functions. Assume that $A_{f,p}$ is locally smaller than $A_{g,q}$. Then
\Eq{1NC}{
  \frac{p_i}{p_0}=\frac{q_i}{q_0} \qquad(i\in\{1,\dots,n\}).
}
If, in addition, $f,g$ are twice differentiable and $p,q$ are continuously differentiable, then the function 
\Eq{fgpq}{
  \frac{q_0^2|g'|}{p_0^2|f'|}
}
is increasing. \\\indent On the other hand, if $f,g$ and $p,q$ are twice continuously differentiable, \eq{1NC} holds and the function \eq{fgpq} has a positive derivative, then $A_{f,p}$ is locally smaller than $A_{g,q}$.}

\begin{proof} In view of \lem{DB}, the differentiability of $f,g$ and continuity of $p,q$ implies that the partial derivatives of the means $M:=A_{f,p}$ and $N:=A_{q,q}$ exist at the diagonal points of $I^n$. Therefore, by \thm{PZ17}, we have that \eq{1A} is valid. Thus, using assertion (1) of \lem{DB}, it follows that
\Eq{*}{
  \frac{p_i}{p_0}=\partial_i A_{f,p}^\Delta
  =\partial_i A_{g,q}^\Delta=\frac{q_i}{q_0}\qquad(i\in\{1,\dots,n\}),
}
which shows that \eq{1NC} is valid.
In what follows, let $r_0$ denote the ratio function $q_0/p_0$. Then, according to \eq{1NC}, we have that $q_i=r_0p_i$.

If, in addition, $f,g$ are twice differentiable and $p,q$ are continuously differentiable, then \lem{DB} implies that the second-order partial derivatives of the means $M:=A_{f,p}$ and $N:=A_{q,q}$ exist at the diagonal points of $I^n$. Therefore, by \thm{PZ17}, we have that the matrix given by \eq{1B} is positive semidefinite. This yields that
\Eq{2NC}{
\partial_i^2 (A_{g,q}-A_{f,p})^\Delta\geq0.
}
Applying assertion (2a) of \lem{DB} and the equality $q_i=r_0p_i$, for $i\in\{1,\dots,n\}$, we have that
\Eq{*}{
  \partial_i^2 (A_{g,q}-A_{f,p})^\Delta&=
  2\frac{q_i'(q_0-q_i)}{q_0^2}+\frac{q_i(q_0-p_i)}{q_0^2}\cdot\frac{g''}{g'}
  -2\frac{p_i'(p_0-p_i)}{p_0^2}-\frac{p_i(p_0-p_i)}{p_0^2}\cdot\frac{f''}{f'}\\
  &=2\frac{(r_0'p_i+r_0p_i')(p_0-p_i)}{r_0p_0^2}+\frac{p_i(p_0-p_i)}{p_0^2}\cdot\frac{g''}{g'}
  -2\frac{p_i'(p_0-p_i)}{p_0^2}-\frac{p_i(p_0-p_i)}{p_0^2}\cdot\frac{f''}{f'}\\
  &=\frac{p_i(p_0-p_i)}{p_0^2}\bigg(2\frac{r_0'}{r_0}+\frac{g''}{g'}-\frac{f''}{f'}\bigg)
  .
}
The factor $\frac{p_i(p_0-p_i)}{p_0^2}$ is positive, therefore, the inequality \eq{2NC} simplifies to 
\Eq{NN}{
  0\leq2\frac{r_0'}{r_0}+\frac{g''}{g'}-\frac{f''}{f'}
  =\bigg(\log\bigg(\frac{q_0^2|g'|}{p_0^2|f'|}\bigg)\bigg)'.
}
This proves that the function \eq{fgpq} is increasing.

In the last part of the proof, assume that $f,g$ and $p,q$ are twice continuously differentiable, \eq{1NC} holds and the function \eq{fgpq} has a positive derivative. Then the means $M:=A_{f,p}$ and $N:=A_{q,q}$ are twice continuously differentiable on $I^n$ and \eq{NN} is valid with strict inequality everywhere on $I$.
We show that the matrix \eq{1B} is positive definite. The diagonal entries of this matrix have been computed above. For the entries off the diagonal, i.e., for $i,j\in\{1,\dots,n\}$ with $i\neq j$, according assertion (2b) of \lem{DB} and the equalities $q_i=r_0p_i$ and and $q_j=r_0p_j$, we have
\Eq{*}{
  \partial_i\partial_j (A_{g,q}-A_{f,p})^\Delta
  &=-\frac{(q_iq_j)'}{q_0^2}-\frac{q_iq_j}{q_0^2}\cdot\frac{g''}{g'}
  +\frac{(p_ip_j)'}{p_0^2}+\frac{p_ip_j}{p_0^2}\cdot\frac{f''}{f'}\\
  &=-\frac{2r_0r_0'p_ip_j+r_0^2(p_ip_j)'}{r_0^2p_0^2}-\frac{p_ip_j}{p_0^2}\cdot\frac{g''}{g'}
  +\frac{(p_ip_j)'}{p_0^2}+\frac{p_ip_j}{p_0^2}\cdot\frac{f''}{f'}\\
  &=-\frac{p_ip_j}{p_0^2}\bigg(2\frac{r_0'}{r_0}+\frac{g''}{g'}-\frac{f''}{f'}\bigg).
}
Therefore, for $i,j\in\{1,\dots,n\}$, 
\Eq{*}{
  \partial_i\partial_j (A_{g,q}-A_{f,p})^\Delta
  =\frac{p_i(\delta_{i,j}p_0-p_j)}{p_0^2}\bigg(2\frac{r_0'}{r_0}+\frac{g''}{g'}-\frac{f''}{f'}\bigg).
}
(Here $\delta_{i,j}$ denotes to Kronecker symbol.)
In order to prove that the matrix \eq{1B} is positive definite, according to Sylvester's criterion, it is necessary and sufficient that the $k\times k$-minors of this matrix should have a positive determinant. 

For $k\in \{1,\dots,n\}$, we obtain
\Eq{*}{
  \det\Big(\big(\partial_i\partial_j (A_{g,q}-A_{f,p})^\Delta\big)_{i,j=1}^k\Big)
  &=\det\Bigg(\bigg(\frac{p_i(\delta_{i,j}p_0-p_j)}{p_0^2}\bigg(2\frac{r_0'}{r_0}+\frac{g''}{g'}-\frac{f''}{f'}\bigg)\bigg)_{i,j=1}^k\Bigg)\\&=\bigg(2\frac{r_0'}{r_0}+\frac{g''}{g'}-\frac{f''}{f'}\bigg)^k\frac{(p_1\cdots p_k)^2}{p_0^{2k}}\det\Bigg(\bigg(\frac{\delta_{i,j}p_0}{p_j}-1\bigg)_{i,j=1}^k\Bigg).
}
Now we compute the last determinant of the above expression:
\Eq{*}{
  \det&\Bigg(\bigg(\frac{\delta_{i,j}p_0}{p_j}-1\bigg)_{i,j=1}^k\Bigg)
  =\begin{vmatrix}
 \dfrac{p_0}{p_1}-1 & -1 &\dots &-1 &-1\\
 -1 &\dfrac{p_0}{p_2}-1 &\dots &-1 &-1\\
  \vdots & \vdots & \ddots& \vdots &\vdots \\
 -1 & -1 & \dots & \dfrac{p_0}{p_{k-1}}-1& -1\\
 -1 & -1 & \dots & -1& \dfrac{p_0}{p_k}-1\\
\end{vmatrix}\\
  &=\begin{vmatrix}
 \dfrac{p_0}{p_1} & 0 &\dots & 0 &-\dfrac{p_0}{p_k}\\
  0 &\dfrac{p_0}{p_2}  &\dots & 0 &-\dfrac{p_0}{p_k}\\
  \vdots & \vdots & \ddots & \vdots &\vdots \\
  0 & 0 & \dots & \dfrac{p_0}{p_{k-1}} & -\dfrac{p_0}{p_k}\\
 -1 & -1 & \dots & -1 & \dfrac{p_0}{p_k}-1\\
\end{vmatrix}
  =\begin{vmatrix}
 \dfrac{p_0}{p_1} & 0 &\dots & 0 &-\dfrac{p_0}{p_k}\\
  0 &\dfrac{p_0}{p_2}  &\dots & 0 &-\dfrac{p_0}{p_k}\\
  \vdots & \vdots & \ddots & \vdots &\vdots \\
  0 & 0 & \dots & \dfrac{p_0}{p_{k-1}} & -\dfrac{p_0}{p_k}\\
 0 & 0& \dots & 0 & \dfrac{p_0}{p_k}\big(1-\big(\dfrac{p_1}{p_0}+\cdots+\dfrac{p_k}{p_0}\big)\big)\\
\end{vmatrix}\\
&=\frac{p_0^k}{p_1\cdots p_k}\bigg(1-\bigg(\frac{p_1}{p_0}+\cdots+\frac{p_k}{p_0}\bigg)\bigg)
=\frac{p_0^{k-1}(p_0-(p_1+\cdots+p_k))}{p_1\cdots p_k}.
}
Therefore,
\Eq{*}{
  \det\Big(\big(\partial_i\partial_j (A_{g,q}-A_{f,p})^\Delta\big)_{i,j=1}^k\Big)
  =\bigg(2\frac{r_0'}{r_0}+\frac{g''}{g'}-\frac{f''}{f'}\bigg)^k\frac{p_1\cdots p_k(p_0-(p_1+\cdots+p_k))}{p_0^{k+1}}.
}
Here, each factor is an everywhere positive function, hence Sylvester's criterion is validated, which then implies that the matrix \eq{1B} is positive definite. In view of the last assertion of \thm{PZ17}, it follows that $M=A_{f,p}$ is locally smaller than $N=A_{q,q}$.
\end{proof}

In what follows, we consider the local comparison problem of nonsymmetric generalized Bajraktarević means when $I$ is an open subinterval of $\R_+$ and the weight functions $p_1,\dots,p_n$ and $q_1,\dots,q_n$ are proportional to power functions.

\Cor{LCp}{Let $I\subseteq\R_+$ and let $f,g:I\to\R$ be differentiable functions with nonvanishing first derivatives. Let $n\geq 2$, $(\lambda_1,\dots,\lambda_n),(\mu_1,\dots,\mu_n)\in\R_+^n$ and $(\alpha_1,\dots,\alpha_n),(\beta_1,\dots,\beta_n)\in\R^n$ and define
\Eq{pq}{ 
p_i(x):=\lambda_i x^{\alpha_i}\qquad\mbox{and}\qquad
q_i(x):=\mu_i x^{\beta_i} 
\qquad(i\in\{1,\dots,n\},\,x\in I).
}
Assume that $A_{f,p}$ is locally smaller than $A_{g,q}$. Then there exist $\gamma>0$ and $\delta\in\R$ such that
\Eq{1NCp}{
  \mu_i=\gamma\lambda_i \qquad\mbox{and}\qquad
  \beta_i=\alpha_i+\delta \qquad(i\in\{1,\dots,n\}).
}
If, in addition, $f$ and $g$ are twice differentiable, then 
\Eq{fgp}{
  x\mapsto x^{2\delta} \frac{|g'(x)|}{|f'(x)|}
}
is increasing on $I$. \\\indent On the other hand, if $f,g$ are twice continuously differentiable, \eq{1NCp} holds and the function \eq{fgp} has a positive derivative, then $A_{f,p}$ is locally smaller than $A_{g,q}$.
}

\begin{proof}
Assume that $f,g:I\to\R$ be differentiable functions with nonvanishing first derivatives and that $A_{f,p}$ is locally nonsmaller than $A_{g,q}$.
Then, by \thm{1NC}, the first-order necessary condition \eq{1NC} holds. This, for all $i,j\in\{1,\dots,n\}$, implies that
\Eq{*}{
  \frac{p_i}{p_j}=\frac{q_i}{q_j},
}
that is,
\Eq{*}{
  \frac{\lambda_i x^{\alpha_i}}{\lambda_j x^{\alpha_j}}
  =\frac{\mu_i x^{\beta_i}}{\mu_j x^{\beta_j}} \qquad(x\in I).
}
Therefore,
\Eq{*}{
  x^{\alpha_i-\alpha_j+\beta_j-\beta_i}
  =\frac{\mu_i\lambda_j}{\mu_j\lambda_i} \qquad(x\in I).
}
The left hand side is a power function which is constant on $I$, hence $\beta_i-\alpha_i=\beta_j-\alpha_j$ and $\mu_i/\lambda_i=\mu_j/\lambda_j$. Hence, there exist $\gamma>0$ and $\delta\in\R$ such that \eq{1NCp} holds.

If, additionally, $f,g:I\to\R$ are twice differentiable functions, then by the second-order necessary condition of \thm{1NC}, the function \eq{fgpq} is increasing. On the other hand, using \eq{1NC} and then \eq{1NCp}, we have
\Eq{id}{
   \frac{q_0^2|g'|}{p_0^2|f'|}(x)
   =\frac{q_i^2|g'|}{p_i^2|f'|}(x)
   =\frac{\mu_i^2x^{2\beta_i}|g'(x)|}
      {\lambda_i^2x^{2\alpha_i}|f'(x)|}
   =\gamma^2 x^{2\delta} \frac{|g'(x)|}{|f'(x)|}.
}
Therefore, the function \eq{fgp} must be increasing.

To prove the reversed statement, suppose that $f,g$ are twice continuously differentiable, \eq{1NCp} holds and the function \eq{fgp} has a positive derivative on $I$. 

Then, for $i\in\{1,\dots,n\}$ and $x\in I$,
\Eq{*}{
  \frac{q_i}{q_0}(x)
  =\frac{\mu_ix^{\beta_i}}{\mu_1x^{\beta_1}+\dots+\mu_nx^{\beta_n}}
  =\frac{\gamma\lambda_ix^{\delta+\alpha_i}}{\gamma\lambda_1x^{\delta+\alpha_1}+\dots+\gamma\lambda_n x^{\delta+\alpha_n}}
  =\frac{\lambda_ix^{\alpha_i}}{\lambda_1x^{\alpha_1}+\dots+\lambda_n x^{\alpha_n}}
  =\frac{p_i}{p_0}(x),
}
which shows that \eq{1NC} is valid on $I$. Furthermore, applying the identity \eq{id}, we can see that the function \eq{fgpq} has a positive derivative. Therefore, in view of the second part of \thm{1NC}, it follows that $A_{f,p}$ is locally smaller than $A_{g,q}$.
\end{proof}

As an immediate consequence of the above corollary, we can characterize the local comparison of nonsymmetric generalized power means.

\Cor{LCpp}{Let $I\subseteq\R_+$ and $a,b\in\R$. Define $f,g:I\to\R$ by
\Eq{fg}{
  f(x):=\begin{cases} 
        x^a & \mbox{if }a\neq0,\\[2mm]
        \log(x) & \mbox{if }a=0,\\
        \end{cases}
        \qquad\mbox{and}\qquad
  g(x):=\begin{cases} 
        x^b & \mbox{if }b\neq0,\\[2mm]
        \log(x) & \mbox{if }b=0.\\
        \end{cases}
} 
Let $n\geq 2$, $(\lambda_1,\dots,\lambda_n),(\mu_1,\dots,\mu_n)\in\R_+^n$, $(\alpha_1,\dots,\alpha_n),(\beta_1,\dots,\beta_n)\in\R^n$ and define $p,q:I\to\R^n_+$ by \eq{pq}. If $A_{f,p}$ is locally smaller than $A_{f,p}$, then there exist $\gamma>0$ and $\delta\in\R$ such that \eq{1NCp} holds and $a\leq b+2\delta$. 
\\ \indent On the other hand, if \eq{1NCp} holds and $a<b+ 2\delta$, then $A_{f,p}$ is locally smaller than $A_{g,q}$.}

\begin{proof} Obviously, \cor{LCp} can be applied and now the functions $f$ and $g$ are infinitely many times differentiable with nonvanishing first derivatives. To show the necessity of the conditions, observe that, for some positive constant $c$,
\Eq{*}{
  x^{2\delta} \frac{|g'(x)|}{|f'(x)|}
  =cx^{b-a+2\delta} \qquad(x\in I).
}
Therefore, \eq{fgp} is increasing if and only if $a\leq b+2\delta$. For the sufficiency, observe that the function \eq{fgp} has a positive derivative if and only if $a< b+2\delta$.
\end{proof}

\section{Global comparison of Bajraktarević means}

If $A_{f,p}$ is globally smaller than $A_{g,q}$, then $A_{f,p}$ is locally smaller than $A_{g,q}$ and hence, in view of \thm{1NC}, the necessity of \eq{1NC} follows and the function \eq{fgpq} must be increasing. However, these conditions are not strong enough to imply the global comparability of $A_{f,p}$ and $A_{g,q}$.

\Thm{2}{Let $f,g:I\to\R$ be strictly monotone differentiable functions with nonvanishing first derivatives and $p,q:I\to\R_+^n$. Assume that \eq{1NC} holds and
\Eq{GSC}{
  \frac{p_0(x)(f(x)-f(y))}{p_0(y)f'(y)}
  \leq \frac{q_0(x)(g(x)-g(y))}{q_0(y)g'(y)} \qquad(x,y\in I).
}
Then $A_{f,p}$ is globally smaller than $A_{g,q}$.}

\begin{proof} Multiplying both sides of \eq{GSC} by $\frac{p_i(x)}{p_0(x)}=\frac{q_i(x)}{q_0(x)}$ (which is valid according to \eq{1NC}), it follows that 
\Eq{GSCi}{
  \frac{p_i(x)(f(x)-f(y))}{p_0(y)f'(y)}
  \leq \frac{q_i(x)(g(x)-g(y))}{q_0(y)g'(y)} \qquad(x,y\in I,\,i\in\{1,\dots,n\}).
}

To show that \eq{comp} holds, let $x_1,\dots,x_n,y\in I$ be arbitrary. Applying the $i$th inequality in \eq{GSCi} at $(x_i,y)$, then adding up the inequalities so obtained side by side, we get
\Eq{GSC1n}{
  \sum_{i=1}^n\frac{p_i(x_i)(f(x_i)-f(y))}{p_0(y)f'(y)}
  \leq \sum_{i=1}^n\frac{q_i(x_i)(g(x_i)-g(y))}{q_0(y)g'(y)}.
}
Let us now take the substitution $y:=A_{f,p}(x_1,\dots,x_n)$ in the above inequality. Then, by the definition of the mean $A_{f,p}$, we have
\Eq{*}{
  f(y)=\frac{\sum_{i=1}^n p_i(x_i)f(x_i)}{\sum_{i=1}^n p_i(x_i)}.
}
Therefore,
\Eq{*}{
  \sum_{i=1}^n\frac{p_i(x_i)(f(x_i)-f(y))}{p_0(y)f'(y)}
  &=\frac{1}{p_0(y)f'(y)}\bigg(\sum_{i=1}^n p_i(x_i)f(x_i)
  -f(y)\sum_{i=1}^n p_i(x_i)\bigg)\\
  &=\frac{1}{p_0(y)f'(y)}\bigg(\sum_{i=1}^n p_i(x_i)f(x_i)
  -\frac{\sum_{i=1}^n p_i(x_i)f(x_i)}{\sum_{i=1}^n p_i(x_i)}\sum_{i=1}^n p_i(x_i)\bigg)=0.
}
Thus, the inequality \eq{GSC1n} simplifies to
\Eq{*}{
  0\leq \sum_{i=1}^n\frac{q_i(x_i)(g(x_i)-g(y))}{q_0(y)g'(y)}.
}
Assume first that $g'$ is positive everywhere in $I$. Then $g$ and also its inverse are strictly increasing and the previous inequality is equivalent to
\Eq{*}{
  0\leq \sum_{i=1}^n q_i(x_i)(g(x_i)-g(y)),
}
which implies that
\Eq{*}{
  g(y)\sum_{i=1}^n q_i(x_i)\leq \sum_{i=1}^n q_i(x_i)g(x_i).
}
Dividing this inequality by $\sum_{i=1}^n q_i(x_i)$ and then applying $g^{-1}$ to the inequality side by side, it follows that $y\leq A_{q,q}(x_1,\dots,x_n)$,
that is, the inequality \eq{comp} is fulfilled. 

If $g'$ is everywhere negative, then $g$ and also its inverse are strictly decreasing and the proof is similar with obvious modifications.
\end{proof}

\Rem{1}{In this remark we show that the sufficient condition \eq{GSC} implies that the function \eq{fgpq} is increasing. We may assume that $f'$ and $g'$ are positive (and hence $f$ and $g$ are strictly increasing). Let $x,y\in I$ be arbitrary with $x<y$. Then, \eq{GSC} implies that
\Eq{xy}{
    \frac{p_0(x)q_0(y)g'(y)}{p_0(y)q_0(x)f'(y)}
  \geq \frac{g(x)-g(y)}{f(x)-f(y)}.
}
On the other hand, interchanging the roles of $x$ and $y$ in \eq{GSC}, we get
\Eq{*}{
    \frac{p_0(y)(f(y)-f(x))}{p_0(x)f'(x)}
  \leq \frac{q_0(y)(g(y)-g(x))}{q_0(x)g'(x)}.
}
This implies that
\Eq{yx}{
    \frac{p_0(y)q_0(x)g'(x)}{p_0(x)q_0(y)f'(x)}
  \leq \frac{g(y)-g(x)}{f(y)-f(x)}.
}
Combining the inequalities \eq{xy} and \eq{yx}, it follows that
\Eq{*}{
  \frac{p_0(y)q_0(x)g'(x)}{p_0(x)q_0(y)f'(x)}
  \leq \frac{g(x)-g(y)}{f(x)-f(y)}
  \leq \frac{p_0(x)q_0(y)g'(y)}{p_0(y)q_0(x)f'(y)},
}
and hence
\Eq{*}{
  \frac{q_0(x)^2g'(x)}{p_0(x)^2f'(x)}
  \leq \frac{q_0(y)^2g'(y)}{p_0(y)^2f'(y)},
}
which proves that the function \eq{fgpq} is increasing. The proof in the cases when at least one of the functions $f'$ and $g'$ is negative is completely analogous.}

In what follows, we present two particular cases when \eq{GSC} is equivalent to the increasingness of the function \eq{fgpq}. In the first setting, the weight functions $p$ and $q$ coincide.

\Thm{2+}{Let $f,g:I\to\R$ be strictly monotone differentiable functions with nonvanishing first derivatives and $p,q:I\to\R_+^n$. Assume that \eq{1NC} holds and the functions
\Eq{GSC+}{
  \frac{q_0}{p_0} \qquad\mbox{and}\qquad \frac{|g'|}{|f'|}
}
are increasing. Then $A_{f,p}$ is globally smaller than $A_{g,q}$.}

\begin{proof} To verify this statement, it suffices to show that the increasingness of the functions \eq{GSC+} implies the inequality \eq{GSC}. We assume that $f'$ and $g'$ are positive (and hence $f$ and $g$ are strictly increasing).

Let $x,y\in I$ be arbitrary. If $x=y$, then \eq{GSC} is trivial. We consider now the case when $x<y$. Then, by the increasingness of the first function in \eq{GSC+}, we have
\Eq{*}{
  \frac{q_0(x)}{p_0(x)}\leq \frac{q_0(y)}{p_0(y)}.
}
By the Cauchy Mean Value Theorem, there exists $z\in\,]x,y[\,$ such that
\Eq{*}{
  0<\frac{g(x)-g(y)}{f(x)-f(y)}
  =\frac{g'(z)}{f'(z)}\leq \frac{g'(y)}{f'(y)},
}
where, for the inequality, we used that $z<y$ and $g'/f'$ is increasing. Multiplying the respective sides of the above inequalities by $q_0(x)/p_0(x)$, it follows that
\Eq{*}{
  \frac{q_0(x)(g(x)-g(y))}{p_0(x)(f(x)-f(y))}
  \leq \frac{q_0(y)g'(y)}{p_0(y)f'(y)}.
}
This inequality, using that $f(x)-f(y)<0$ implies that \eq{GSC} is satisfied if $x<y$. In the case when $y<x$, a completely analogous argument shows that the inequality \eq{GSC} is also true.

The proof when $f'$ or $g'$ is negative is completely similar.
\end{proof}

\Thm{3}{Let $f,g:I\to\R$ be strictly monotone twice continuously differentiable functions with nonvanishing first derivatives and $p:I\to\R_+^n$ be a continuous function. Then the following conditions are pairwise equivalent.
\begin{enumerate}[(i)]
 \item $A_{f,p}$ is globally smaller than $A_{g,p}$;
 \item $A_{f,p}$ is locally smaller than $A_{g,p}$;
 \item The function $|g'/f'|$ is increasing on $I$;
 \item \Eq{*}{\frac{f''}{f'}\leq \frac{g''}{g'};}
 \item Provided that $g$ is increasing (decreasing), the function $g\circ f^{-1}$ is convex (concave) on $f(I)$;
 \item \Eq{GSCfg}{
  \frac{f(x)-f(y)}{f'(y)}
  \leq \frac{g(x)-g(y)}{g'(y)} \qquad(x,y\in I).}
\end{enumerate}}

\begin{proof}
 The implication (i)$\Rightarrow$(ii) is trivial. Using \thm{1NC}, (iii) is an immediate consequence of (ii). 
 
 Now assume that (iii) is valid, i.e., $|g'/f'|$ is increasing on $I$. In the case $g'/f'>0$ this implies that the function $(g'/f')'=(g'/f')(g''/g'-f''/f')$ is nonnegative, whence (iv) follows. In the case $g'/f'<0$ an analogous argument yields that (iv) is also true.
 
 Suppose that (iv) holds and $g$ is increasing (the other case can be treated very similarly). Since  
 \Eq{2der}{
	(g\circ f^{-1})''=((g\circ f^{-1})')'=\bigg(\frac{g'}{f'}\circ f^{-1}\bigg)'
	= \bigg(\frac{g'}{f'^2}\cdot\bigg(\frac{g''}{g'}-\frac{f''}{f'}\bigg)\bigg)\circ f^{-1}\geq0,
}
the function $g\circ f^{-1}$ is convex.

Assume now that (v) holds and $g$ is increasing (the other case is completely similar). Then, by the convexity of $h:=g\circ f^{-1}$, for all $u,v\in f(I)$, we have $h(v)+h'(v)\cdot(u-v)\leq h(u)$, that is,
\Eq{*}{
 g\circ f^{-1}(v)+\bigg(\frac{g'}{f'}\circ f^{-1}\bigg)(v)\cdot(u-v)\leq g\circ f^{-1}(u).
}
Substituting $u:=f(x)$ and $v:=f(y)$, where $x,y\in I$, it follows that
\Eq{*}{
 g(y)+\frac{g'(y)}{f'(y)}(f(x)-f(y))\leq g(x),
}
which shows that (vi) is satisfied. 

Finally, suppose that (vi) is valid. Then, the inequality \eq{GSC} is also true with $q_0:=p_0$. Therefore, according to \thm{2}, we obtain that $A_{f,p}$ is globally smaller than $A_{g,p}$, i.e., (i) is fulfilled. 
\end{proof}

In the second setting, the functions $f$ and $g$ are the same.

\Thm{4}{Let $f:I\to\R$ be strictly monotone twice continuously differentiable function with a nonvanishing first derivative and $p,q:I\to\R_+^n$ be a continuous functions. Then the following conditions are pairwise equivalent.
\begin{enumerate}[(i)]
 \item $A_{f,p}$ is globally smaller than $A_{f,q}$;
 \item $A_{f,p}$ is locally smaller than $A_{f,q}$;
 \item The condition \eq{1NC} holds and the function $q_0/p_0$ is increasing on $I$.
\end{enumerate}}

\begin{proof}
The implication (i)$\Rightarrow$(ii) is obvious. 

Assume first that (ii) is valid. Then, according to \thm{1NC}, the equalities in \eq{1NC} must be satisfied and the function $(q_0^2|g'|)/(p_0^2|f'|)=q_0^2/p_0^2$ is increasing. Therefore, condition (iii) holds true.

Finally, suppose that (iii) is satisfied. Then, with $g:=f$, we can see that the functions in \eq{GSC+} are increasing.
Hence, according to \thm{2+}, $A_{f,p}$ is globally smaller than $A_{g,q}$ and assertion (i) is valid.
\end{proof}

\Cor{2}{Let $I\subseteq\R_+$ and let $f,g:I\to\R$ be differentiable functions with nonvanishing first derivatives. Let $n\geq 2$, $(\lambda_1,\dots,\lambda_n),(\mu_1,\dots,\mu_n)\in\R_+^n$ and $(\alpha_1,\dots,\alpha_n),(\beta_1,\dots,\beta_n)\in\R^n$ and define $p,q:I\to\R_+^n$ by \eq{pq}. Assume that there exist $\gamma>0$ and $\delta\in\R$ such that \eq{1NCp} is satisfied and
\Eq{GSCp}{
  \frac{f(x)-f(y)}{f'(y)}
  \leq \frac{x^\delta(g(x)-g(y))}{y^\delta g'(y)} \qquad(x,y\in I).
}
Then $A_{f,p}$ is globally smaller than $A_{g,q}$.}

\begin{proof} Due to the definition of $p$ and $q$ by \eq{pq} and conditions \eq{1NCp}, we have that \eq{1NC} is valid. We are going to prove that the inequality \eq{GSC} is also satisfied. First, observe that \eq{pq} and conditions \eq{1NCp} imply that $q_0(x)=\gamma p_0(x)x^\delta$ for all $x\in I$. Therefore, multiplying the inequality \eq{GSCp} by $p_0(x)/p_0(y)$ side by side, we get
\Eq{*}{
  \frac{p_0(x)(f(x)-f(y))}{p_0(y)f'(y)}
  \leq \frac{\gamma p_0(x)x^\delta(g(x)-g(y))}{\gamma p_0(y)y^\delta g'(y)}
  = \frac{q_0(x)(g(x)-g(y))}{q_0(y)g'(y)},
}
which shows that the inequality \eq{GSC} is satisfied. Hence, using \thm{2}, it follows that  
$A_{f,p}$ is globally smaller than $A_{g,q}$. 
\end{proof}

\Cor{3}{Let $I\subseteq\R_+$, $n\geq 2$, $(\lambda_1,\dots,\lambda_n),(\mu_1,\dots,\mu_n)\in\R_+^n$, $(\alpha_1,\dots,\alpha_n),(\beta_1,\dots,\beta_n)\in\R^n$ and $a,b\in\R$. Define $p,q:I\to\R_+^n$ by \eq{pq} and $f,g:I\to\R$ by \eq{fg}. Assume that there exist $\gamma>0$ and $\delta\in\R$ such that \eq{1NCp} is satisfied and
\Eq{GSCpp}{
  \min(a,0)\leq \delta+\min(b,0) \qquad\mbox{and}\qquad
  \max(a,0)\leq \delta+\max(b,0).
}
Then $A_{f,p}$ is globally smaller than $A_{g,q}$.}

\begin{proof} Due to the definition of $p$ and $q$ by \eq{pq} and conditions \eq{1NCp}, we have that \eq{1NC} is valid. We are going to prove that the inequality \eq{GSCpp} is also satisfied.

For $r,s\in\R$, define
\Eq{*}{
  g_{r,s}(t)
  :=\begin{cases}
  \dfrac{t^r-t^s}{r-s}&\mbox{if } r\neq s,\\[2mm]
  t^r\log(t)&\mbox{if } r=s.
  \end{cases}
}
The mapping $r\mapsto t^r$ is convex for all $t>0$. Therefore, if 
\Eq{*}{
  \min(r,s)\leq\min(u,v) \qquad\mbox{and}\qquad
  \max(r,s)\leq\max(u,v),
}
then
\Eq{*}{
  g_{r,s}(t)\leq g_{u,v}(t) \qquad(t\in\R_+). 
}
Thus, the inequalities in \eq{GSCpp} imply that
\Eq{ab}{
  g_{a,0}(t)\leq g_{b+\delta,\delta}(t) \qquad(t\in\R_+).
}
If $ab\neq0$, this inequality is equivalent to
\Eq{*}{
  \dfrac{t^a-1}{a}\leq t^\delta\dfrac{t^b-1}{b} \qquad(t\in\R_+).
}
For $x,y\in I$, substitute $t:=x/y$ into this inequality.
After simplifications, we obtain
\Eq{*}{
  \dfrac{x^a-y^a}{ay^{a-1}}\leq \Big(\frac{x}{y}\Big)^\delta\dfrac{x^b-y^b}{by^{b-1}}.
}
This inequality, according to the definition \eq{fg} of $f$ and $g$, can be rewritten as \eq{GSCpp}. In the case when $ab=0$, a similar argument shows that \eq{ab} also implies \eq{GSCpp}. 

Now, we are able to apply \cor{2}, and hence we can conclude the result.
\end{proof}



\providecommand{\bysame}{\leavevmode\hbox to3em{\hrulefill}\thinspace}
\providecommand{\MR}{\relax\ifhmode\unskip\space\fi MR }
\providecommand{\MRhref}[2]{%
  \href{http://www.ams.org/mathscinet-getitem?mr=#1}{#2}
}
\providecommand{\href}[2]{#2}

\end{document}